\input amstex.tex
\documentstyle{amsppt}
\document
\NoBlackBoxes
\NoRunningHeads


             \def \lra {\longrightarrow}
             
             \def \bD {{\Bbb D}}

\def\end{\text{End}}

\def\spf{\text{Spf}}
\define\bH{\Bbb H}
             
             \def\wt{\text{wt}\,}
             \def\tw{\hbox{\rm Tw}}
             \def\tws{{\buildrel\sim\over{\hbox{\rm Tw}}}}

             \def \cH {{\Cal H}}
             \def \cA {{\Cal A}}
             
             \def \cO {{\Cal O}}
             \def \OK {{\Cal O}_K}
             
             \def \cS {{\Cal S}}
             \def \cW {{\Cal W}}
             \def \cX {{\Cal X}}

             \def \Ro {{R_{ord}^\prime}}
             \def \vpo {\varphi_{ord}^\prime}
             \def \cXoa {{\Cal X}_{ord}}
             \def \Xoa {X_{ord}}
             \def \Roa {{R_{ord}}}
             \def \vpoa {\varphi_{ord}}
             
             \def \cXoe {\cXoa^{esd}}
             \def \X {X}
             \def \Xp {X_0}

             \def \oQ {{\overline\Q}}
\def \OoQ {{\Cal O}_{\oQ_p}}


\define\z{\Bbb Z}
\define\A{\Bbb A}

\define\I{\Bbb I}
\define\R{\Bbb R}
\define\C{\Bbb C}
\define\Q{\Bbb Q}
\define\Z{\Bbb Z}

\define\dee{\Bbb D}

\define\inv{^{-1}}

\def\Hom{\text{Hom}}

\def\diag{\text{diag}\,}

\def\ker{\text{ker}}

\def\mod{\text{mod }}
\def\det{\text{det}}

\topmatter
\title{Rigidity of $p$-adic cohomology classes\\
of congruence subgroups of $GL(n,\Bbb Z)$}
\endtitle

\author{Avner Ash*, David Pollack**, and Glenn Stevens***}
\endauthor

\affil{*Department of Mathematics, Boston College, Chestnut Hill MA
02467
\smallskip
**Department of Mathematics, Wesleyan University, Middletown CT 06459
\smallskip
***Department of Mathematics, Boston University, Boston, MA 02215}
\endaffil
\thanks{First author partially supported by NSA grant MDA 904-00-1-0046
and NSF grant
DMS-0139287. Third author partially supported by NSF grant DMS 0071065.
This manuscript is submitted for publication with the understanding that
the United States government is authorized to reproduce and distribute
reprints.}
\endthanks

\endtopmatter

The main result of [Ash-Stevens 97] describes a framework for
constructing
$p$-adic analytic families of $p$-ordinary arithmetic Hecke
eigenclasses in the
cohomology of congruence subgroups of $GL(n)/\Q$ where the Hecke
eigenvalues
vary $p$-adic analytically as functions of the weight. Unanswered in
that paper
was the question of whether or not {\it every}
$p$-ordinary arithmetic eigenclass can be ``deformed" in a ``positive
dimensional family" of arithmetic eigenclasses. In this paper we make
precise
the notions of ``deformation" and ``rigidity" and investigate their
properties.
Rigidity corresponds to the non-existence of positive dimensional
deformations
other than those coming from twisting by the powers of the determinant.
Formal
definitions are given in $\S$3. When $n=3$, we will give a necessary and
sufficient condition for a given $p$-ordinary arithmetic eigenclass to
be
$p$-adically rigid and we will use this criterion to give examples of
$p$-adically rigid eigenclasses on $GL(3)$.

 On the other hand, an argument of Hida in section 5 of [Hida 94] 
can be modified to prove the existence of
nontrivial, non-$p$-torsion $p$-adic deformations (modulo twisting) for
the same examples.  Details will appear elsewhere.  Here we merely
note that the use of quasicuspidality (see Definition 6.1) seems to be
essential to make Hida's arguments work in the $GL(3)$ case.

 Coupled with our rigidity result (Theorem 9.4) this proves the
 existence of non-arithmetic $p$-adic deformations of $p$-ordinary
 arithmetic (cohomological) cuspforms for $GL(3)$.  In fact, both of
our
 examples are deformations of Hecke eigenclasses in the cohomology of
 congruence subgroups of $SL(3,\Z)$ with trivial coefficients.

More precisely, we define a ``Hecke eigenpacket" to be a map from the
Hecke
algebra to a  coefficient ring which gives the Hecke eigenvalues
attached to a
Hecke eigenclass in the cohomology.  We then consider $p$-adic analytic
``deformations" of the Hecke eigenpacket.  There are several advantages
to this point of view, as opposed to deforming the cohomology classes
themselves.  For one thing, a surjective map of coefficient modules need
not induce a surjective map on the cohomology, but the long exact 
sequence in cohomology enables us to keep track of Hecke eigenpackets.  
For example, Lemma 4.3 is a surjectivity result on the level of Hecke 
eigenpackets proved in this way (see Theorem 5.1).  Standard tools of 
algebra
and geometry can be used to study the Hecke algebra and the deformations
of the eigenpackets.  Moreover, it is the Hecke eigenvalues that carry
the interesting arithmetic information attached to a cohomological
eigenclass.

Our main theoretical tool is the cohomology of congruence subgroups
$\Gamma\subseteq GL(n,\z)$ with coefficients in the space of bounded
measures
on the big cell of $GL(n)$ (see $\S$2). In $\S$4 we define the ring
$\Roa$ as a
quotient of the Hecke algebra acting on the $p$-ordinary part of this
cohomology. In Theorem 4.2 we prove that $\Roa$ is {\it universal} in
the sense
that every $p$-ordinary arithmetic Hecke eigenpacket factors through
$\Roa$. We
show that $\Roa$ is a complete semilocal noetherian ring and study the
geometry
of its associated
$p$-adic rigid analytic space $\cXoa$. This leads us to a
necessary and sufficient condition for $p$-adic rigidity (Theorem 4.9),
valid
for any $n$. In sections 6 through 8 we recall a number of
representation
theoretic facts, which in case $n=3$ allow us to prove a criterion for
rigidity
(Corollary 8.5) that is accessible to numerical computation. Our
particular
computations are described in $\S$9.

The first and third authors have also developed a theory of $p$-adic
deformations for non-ordinary (but finite slope) Hecke eigenpackets.
Using a
more representation-theoretic approach, [Emerton 04] has defined an
eigenvariety in this context, generalizing to $GL(n)$ the Coleman-Mazur
eigencurve for $GL(2)$. 
Emerton states his results for any $\Q_p$-split reductive
algebraic group and it is perhaps worth noting that our methods also
apply {\it
mutatis mutandis} in this more general context.

The question of
$p$-adic rigidity arises again in the finite slope context. However, it
will probably
be impossible to verify by computations like those in this paper
whether or not any
particular non-ordinary eigenpacket is
$p$-adically rigid. The reason is that in the non-ordinary case, we
lose control over
the ``size" of the subset of the parameter (``weight") space over which
the family
exists. In the ordinary case, we know that the families exist over the
whole
parameter space, so we can test
$p$-adic rigidity by looking at weights that are
$p$-adically far from the given weight. Since the only weights for
which computer
calculations are currently feasible are $p$-adically rather far from
each other, this
limits our computational methods to the ordinary case.

One knows that the symmetric square lift to $GL(3)$ of any
$p$-ordinary classical newform $f$ does indeed deform into a $p$-adic
analytic
family (lifting Hida's $p$-ordinary deformation of $f$). From
[Ash-Stevens 86] we know
that the cohomology classes associated to the symmetric square lifts are
essentially self-dual eigenclasses (Definition 8.1) for $GL(3)$. In
looking for examples
of $p$-adically rigid eigenclasses we therefore examined the
non-essentially-self-dual
$p$-ordinary eigenclasses constructed numerically in [Ash-Grayson-Green
84]. The
numerical calculations described in $\S$9 show that these eigenclasses
satisfy our
criterion and are therefore $p$-adically rigid.

The calculations suggest that, outside of symmetric square lifts,
non-torsion
arithmetic cohomology is very sporadic. On the other hand, failure of
rigidity
of any non-symmetric square lift would give rise to an abundance of
non-essentially-selfdual eigenclasses, at least if the Hecke eigenpacket
is not essentially-selfdual modulo $p$. So we are led to conjecture:

\proclaim{Conjecture 0.1} If a finite slope
cuspidal Hecke eigenclass in the cohomology of a
congruence subgroup of $GL(3,\z)$ is not essentially selfdual then it
is $p$-adically
rigid.
\endproclaim

Recent work of Eric Urban suggests that the converse may be true for
all $GL(n)$, $n\ge
2$. When $n=2$, every cuspidal eigenclass is essentially selfdual, and
in this case
Hida (in the ordinary case) and Coleman (in general) proved that all
finite slope
cuspidal eigenclasses can be deformed into non-trivial $p$-adic
analytic families, i.e
they are not $p$-adically rigid in our sense.

In order to prove $p$-adic rigidity in the cases we study in the last
section of the
paper, we have to compute the cuspidal cohomology of $\Gamma_0(q)$ with
small-dimensional irreducible coefficient modules when $n=3$. Here,
$\Gamma_0(q)$
denotes the arithmetic subgroup of $SL(n,\z)$ consisting of matrices
whose top row is
congruent to $(*,0,\dots,0)$ modulo $q$, where $q$ is a prime different
from $p$.

Here are the computational results.
Let $V_h$ denote the irreducible finite-dimensional $GL(3,\C)$-module of
 highest weight $(2h,h,0)$.  $V_h$ is a submodule of $Sym^h(A) \otimes Sym^h(A^\dagger) \otimes det^h$
where $A$ is the standard 3-dimensional complex representation of
$GL(3,\C)$ and $A^\dagger$ is its contragredient. Let $E$ be a rational
complex
representation of $GL(3,\C)$. As explained in section 8 below, the
cuspidal cohomology
of any congruence subgroup of $SL(3,\z)$ with coefficients in $E$ is
known to vanish
if
$E$ contains no submodule
isomorphic to a twist by a power of the determinant of any $V_h$. So to
look for
cuspidal classes, it suffices to consider cohomology with $V_h$
coefficients.

\proclaim{Computational Results 0.1}

There is no cuspidal cohomology for $\Gamma_0(q)$ with coefficients
in $V_h$ for
$h=1,2$,
$q$ any prime up to and including 223, and for $h=3,4$, $q$ any prime
up to and including 97.

\endproclaim

The computations use the programs and methods outlined in section 8 of
[Ash-Doud-Pollack
02]. Cuspidal cohomology can
only occur in $H^2$ and $H^3$ and Poincar\'e duality shows that
$H_{cusp}^2 \simeq H_{cusp}^3$. So it suffices
to compute $H^3$.

It is somewhat surprising that we found no non-essentially-selfdual
cuspidal cohomology
in these ranges, since for $h=0$ (i.e. trivial coefficients) there
exists
non-essentially-selfdual cuspidal cohomology for $q=53,61,79,89$
[Ash-Grayson-Green 84].

As a consequence of these computations, also in section 9, we deduce
the 3-adic rigidity
of the cuspidal cohomology with trivial coefficients
of level 89 and the 5-adic rigidity of the cuspidal cohomology of level
61. The difficulty of testing
other examples is purely owing to our inability to compute for large
$h$.

In this paper we avoid speaking about an eigenvariety for $GL(n)$,
preferring instead to formulate our results in terms of eigenpackets.
We look only at ``arithmetic eigenpackets'', i.e. those whose
associated rigid analytic spaces are the Zariski closure of their
arithmetic points (see definitions 3.1 and 3.6).  In particular we say
an arithmetic point is ``rigid'' if it is isolated on any {\it
arithmetic}  eigenpacket that contains it (see definition 3.9).

It may help to describe our picture roughly in terms of a putative
eigenvariety.    We focus on the $n=3$ situation and let $Y$ be the
eigenvariety for $GL(3)$ modulo twisting.  $Y$ projects to a
$2$-dimensional weight space $\cW$.  If a point $x$ on $Y$ is a
symmetric square lift we expect it will lie on a component of $Y$
whose projection to $\cW$ contains the $1$-dimensional self-dual line.
If $x$ is not a symmetric square lift, then our conjecture above
implies that any component of $Y$ containing $x$ will project down to a set in 
$\cW$ that is either $1$-dimensional and transverse to the self-dual line or $0$-dimensional.

The contents of this paper are as follows. In Sections 1 and 2 we
introduce the Hecke
algebras and the cohomology groups of interest to us and describe the
action of the Hecke algebra on the cohomology. Section 2 also recalls
the main
results from [Ash-Stevens 97] that form the basis of the rest of this
paper.
In section 3 we develop a formalism of Hecke eigenpackets and define
the key notions
of $p$-adic deformation and $p$-adic rigidity. In section 4, we
introduce the
$p$-ordinary arithmetic eigenvariety and in section 5 we prove its
universality.

In sections 6 through 8 we collect the representation theoretic facts
we need
in order to give a testable criterion for $p$-adic
rigidity. We then give our computational results in Section 9. In
particular, we
give two examples of $p$-adically rigid eigenclasses, one with $p=3$
and level 89, the
other with $p=5$ and level 61.

\heading{Section 1:~~Hecke Operators and $(T,S)$-Modules}
\endheading

Let $G= GL(n)$ ($n\ge 2$) and fix a positive prime $p$. Let $T\subseteq
G$ be the group
of diagonal matrices and let $S\subseteq G(\Q)$ be a subsemigroup
containing a
given congruence subgroup $\Gamma\subseteq SL_n(\z)$.
\medskip

\noindent
{\bf Definition 1.1.} A $(T,S)$-module is a topological $\z_p$-module
$V$
endowed with a continuous left action of $T(\Q_p)$ and a right action
of $S$ that
commutes with the $T(\Q_p)$-action.
\medskip

\noindent
In the examples, $V$ can always be viewed as a ``highest weight module"
in which the action of $T$ is defined via multiplication by the highest
weight
character of $T(\Q_p)$.

Let $\Lambda$ and $\Lambda_T$ be the completed group rings on the
topological groups
$T(\z_p)$ and $T(\Q_p)$ respectively. Then $\Lambda$ is the Iwasawa
algebra of
$T(\z_p)$ and, letting
$T_0\subseteq T(\Q_p)$ be the subgroup of matrices whose diagonal
entries are integral
powers of $p$, we have
$
\Lambda_T = \Lambda[T_0].
$
The group $T_0$ is free abelian with generators
$\pi_k := \diag(1,\ldots,1,p,\ldots,p)$ ($k$ copies of $p$), for
$k=1,\ldots,n$, and consequently we have
$$
\Lambda_T = \Lambda[\pi_1,\pi_1^{-1},\cdots,\pi_n,\pi_n^{-1}],
$$
a Laurent polynomial ring in $n$ variables over $\Lambda$.
The element
$$
\pi := \pi_1\cdots\pi_{n-1}
$$
plays a special role in our theory.

Since $T(\Q_p)$ acts continuously on any $(T,S)$-module $V$, we may
regard
$V$ as a right $\Lambda_T[S]$-module. Indeed, this defines an
equivalence of categories between the category of
$(T,S)$-modules and the category of continuous $\Lambda_T[S]$-modules.

Now fix a positive integer $N$ with $p\not|N$ and define the congruence
subgroups
$$\eqalign{
\Gamma &:= \{g \in SL(n,\z) \mid \text{top row of } g \equiv
(*,0,\dots,0)\
\mod N
\}, \text{and}\cr
\Gamma_0 &:= \{ g\in\Gamma\mid
g~\text{is~upper~triangular~modulo}~p~\}.\cr
}$$
We define the subsemigroups $S^\prime$ and $S_0^\prime$ of $M_n^+(\Bbb
Z)$ (the
plus denotes positive determinant) by
$$\eqalign{
S^\prime &:= \{s \in M_n^+(\z) \mid (\hbox{\rm det} (s), Np) = 1
\text{, top row of } s \equiv (*,0,\dots,0)\ \mod N \}\cr
S'_0 &:= \{ s\in S^\prime\mid s~\hbox{\rm is upper triangular modulo
$p$} \}.\cr
}$$

Finally, we define the subsemigroups $S$ and $S_0$ of $G(\Q_p)$:

$$\eqalign{
S &:=\text{semigroup~generated~by}~S'~\hbox{\rm and $\pi$},\cr
S_0 &:=\text{semigroup~generated~by}~S'_0~\hbox{\rm and $\pi$}.\cr
}$$

For $s\in S$ we let $T(s)$ be the double coset $\Gamma s \Gamma$.
It is well known that the double coset algebra $D_{\Lambda_T}(\Gamma
,S) =
\Lambda_T[\Gamma\backslash S/\Gamma]$ is the free polynomial algebra
over
$\Lambda_T$ generated by $T(\pi)$ and the elements
$$\eqalign{
T(\ell,k) &:= T(s),\ \hbox{\rm $s = \diag (1,\dots ,1,\ell,\dots
,\ell)$,
with
$k$ copies of
   $\ell$'s ($1\le k\le n$)},\cr
}$$
for all primes $l$ not dividing $Np$. (See Props 3.16, 3.30 and Thm 3.20 of [Shimura 71].)

The double coset algebra
$D_{\Lambda_T}(\Gamma_0,S_0)$ is likewise a free polynomial algebra generated by the
$\Gamma_0$-double cosets of $\pi$ and of the
$\diag(1,\dots,1,\ell,\dots,\ell)$ with $l$ a prime not dividing $Np$
and $1\le k\le n$.  Thus the map sending a polynomial expression in
the generators of  $D_{\Lambda_T}(\Gamma,S)$ to the corresponding polynomial
expression in the generators of $D_{\Lambda_T}(\Gamma_0,S_0)$ gives an
isomorphism 
$
D_{\Lambda_T}(\Gamma,S) \lra D_{\Lambda_T}(\Gamma_0,S_0).
$

We define the element $U$ to be $T(\pi)$ 
and let $\cH$ denote the $\Lambda$-subalgebra of $D_{\Lambda_T}(\Gamma
,S)$
generated by the elements $U$ and the $T(\ell,k)$:
$$
\cH = \Lambda\left[U,\ T(\ell,k)\,\biggm|\,\hbox{\rm $\ell\not|Np$ is
prime and
$1\le k \le n$}\,\right].
$$
   We let $\cH_0$ be the image of $\cH$ in
$D_{\Lambda_T}(\Gamma_0,S_0)$
and note that
$$
\cH \buildrel\sim\over \lra \cH_0
$$
is an isomorphism of $\Lambda$-modules.

We define a modified 
action of $S$, 
which we call the
$\star$-action, as
follows. There is a unique multiplicative map $\tau : S\lra T_0$ that
is trivial on
$S^\prime$ and sends $\pi$ to itself.
For any $(T,S)$-module $V$ we define the $\star$-action of
$S$ on $V$ by the formula
$$
w\star \sigma := \tau(\sigma)^{-1}\cdot w \sigma.
$$
for $w\in V$ and $\sigma\in S$. We then define $V^\star$ to be $V$
endowed with
the $\star$-action of $S$ on the right and the given action of
$\Lambda$ on the left.
For any subquotient $L$ of $V^\star$
the double coset algebra
$D_\Lambda(\Gamma,S)$ acts naturally on the $\Gamma$-cohomology of $L$.
We then let $\cH$ act on the cohomology of $L$.

We apply this construction verbatim for any $(T,S_0)$-module $V$ and
similarly obtain
a $\Lambda[S_0]$-module $V^\star$. Then for any
$\Lambda[S_0]$-subquotient $L$
of
$V^\star$ the $\Gamma_0$-cohomology of $L$ receives a natural
$\cH_0$-module structure. Using our fixed isomorphism
$\cH\cong \cH_0$, we may (and do) regard the $\Gamma_0$-cohomology of
$L$ as an
$\cH$-module. For example, if $V$ is a $(T,S)$-module and $L$ is a
subquotient of
$V^\star$ then we may also regard $V$ as a $(T,S_0)$-module and $L$ as
an
$\Lambda[S_0]$-subquotient of $V^\star$. So both $H^\ast(\Gamma,L)$ and
$H^\ast(\Gamma_0,L)$ are endowed with natural $\cH$-module structures.
Note, however, that in general the restriction map
$$
H^\ast(\Gamma,L) \lra H^\ast(\Gamma_0,L)
$$
does {\it not} commute with $U$, though it does commute with
the Hecke operators $T(\ell,k)$ for $\ell\not|Np$.

Note that we use throughout the notation $H^\ast$ to stand for 
$\oplus_k H^k$.

\heading{Section 2:~~$p$-Adic Interpolation of Arithmetic Cohomology}
\endheading

Let $B$ be the Borel
subgroup of upper triangular matrices and denote the group of lower
triangular
matrices by $B^{opp}$.
Let $\lambda$ denote the highest weight, with respect to $(B^{opp},T)$,
of
a finite dimensional irreducible right rational $\Q_p[G]$-module
$V_\lambda$.
Then $V_\lambda$ inherits a natural action of $S$ on the right and we
let
$T(\Q_p)$ act (continuously) on the left via the character $\lambda$.
In a moment, we will also
choose a $\z_p$-lattice $L_\lambda\subseteq V_\lambda$ that is invariant
under $S_0^\prime$ and satisfies
$\pi^{-1} L_\lambda\pi \subseteq L_\lambda$. Thus by the discussion of
the last section,
we obtain an action of
$\cH$ on
$H^\ast(\Gamma,V_\lambda)$, $H^\ast(\Gamma_0,V_\lambda)$, and
$H^\ast(\Gamma_0,L_\lambda)$.

Now let $N^{opp}$ be the
unipotent radical of $B^{opp}$ and let $Y := N^{opp}(\Q_p)\backslash
G(\Q_p)$.
We also let
$$\eqalign{
X &:= \hbox{\rm the image of $G(\z_p)$ in $Y$, and}\cr
\Xp &:= \hbox{\rm the image of $I$ in $Y$,}\cr
}$$
where $I$ is the Iwahori subgroup of $G(\z_p)$ consisting of matrices
that are upper triangular modulo $p$.
The semigroup $S$ acts on $Y$ by multiplication on the right and
$T(\Q_p)$
acts by multiplication on the left.
Both $X$ and $\Xp$ are stable under $T(\Z_p)$. Also $X$ is stable under
the action of
$S^\prime$ and $\Xp$ is stable under $S_0^\prime$.
We define a new action of $S$ on $Y$ by setting for $\sigma\in S$
$$N^{opp}(\Q_p)g\star\sigma=N^{opp}(\Q_p)\tau(\sigma)^{-1}g\sigma,$$
where $\tau\colon S\rightarrow T_0$ is the homomorphism sending $\pi$ to $\pi$ and acting trivally on $S'$ introduced above.
This new $\star$-action is well defined since $\pi$ normalizes $N^{opp}(\Q_p)$.    
Moreover, $X_0$ is stable under the $\star$-action of $S_0$ since every element of $X_0$ is represented by some element $b\in B(\Z_p)$ and $\pi^{-1}B(\Z_p)\pi\subset B(\Z_p)$.
We will henceforth use this $\star$-action when we act $S$ on $Y$ or $S_0$ on $X_0$.

For any open subset $Z\subseteq Y$, we let
$\dee(Z)$ be the space of $\z_p$-valued measures on $Z$, so
$\dee(Z) = \Hom_{\z_p}(C(Z), \z_p)$,
where $C(Z)$ is the space of all compactly supported continuous
$\z_p$-valued
functions on
$Z$.
By functoriality (using the $\star$-action), $\dee(Y)$ inherits a natural $(T,S)$-structure. By
extending
measures by zero, we obtain a $\Lambda[S_0]$-embedding
$\bD(\Xp)\hookrightarrow \bD(Y)$.
In [Ash-Stevens 97], the space $\bD=\bD(X)$ is 
given a $\Lambda[S]$ structure by setting for $\mu\in \bD$ and $f\in C(X)$  
$$(\mu\star s)(f)=\mu((s\star\tilde f)|_{X})$$
where $\tilde f$ is the extension by zero of $f$ to $Y$ and $(s\star\tilde f)(y)=\tilde f(y\star s)$.
Hence the cohomology groups
$
H^\ast(\Gamma_0,\bD(\Xp))
$,
$
H^\ast(\Gamma,\bD)
$
are endowed with natural $\cH$-module structures. Moreover, restriction
of
distributions from $X$ to $\Xp$ induces a $\Lambda[S_0]$-morphism $\bD
\lra \bD(\Xp)$
and thus induces a map
$$
H^\ast(\Gamma,\bD)\lra
H^\ast(\Gamma_0,\bD(\Xp))
$$
that commutes with action of $\cH$.  (To see that this maps commutes with $T(\pi)$ requires an explicit calculation similar to the proof of Proposition 4.2 of [Ash-Stevens 97].)

For any dominant integral weight $\lambda$, we define the map
$$\phi_{\lambda}\colon\bD(Y)\lra V_\lambda$$
by $\phi_{\lambda}(\mu )=\int_{\Xp} v_\lambda x\,d\mu (x)$,
where $v_\lambda$ denotes a fixed highest weight vector in $V_\lambda$.
One checks easily that $\phi_\lambda$ is a morphism of $(T,S)$-modules.
Thus $\phi_\lambda$ induces a $\Lambda[S_0]$-morphism
$$
\phi_\lambda : \bD(\Xp) \lra V_\lambda^\star.
$$
We let $L_\lambda$ be the image of this map and note that since
$\bD(\Xp)$
is compact $L_\lambda$ is a finitely generated $I$-invariant submodule
of
$V_\lambda$, hence is a lattice in $V_\lambda$ (since $V_\lambda$ is
irreducible).
This is the lattice whose existence was asserted in the first paragraph
of this section.
Thus $\phi_\lambda$ induces a natural $\cH$-morphism
$$
\phi_{\lambda}\colon
H^\ast (\Gamma ,\Bbb D) \to H^\ast (\Gamma_0 ,L_\lambda).
$$

For any compact $\Bbb Z_p$-module $H$ with a continuous action of the
operator $U$ on it, we define the ordinary part of $H$ to be
$H^0=\cap U^iH$,
where $i$ runs over all positive integers. Then $H^0$ is the largest
submodule of $H$
on which $U$ acts invertibly, and $H \mapsto H^0$ is an exact functor.

The following theorem is Theorem 5.1 of [Ash-Stevens 97].

\proclaim{Theorem 2.1} $H^\ast (\Gamma,\Bbb D)^0$ is finitely generated
as a $\Lambda$-module. Moreover, for every dominant integral weight
$\lambda$, the
kernel of the map
$$\phi_{\lambda}\colon H^\ast (\Gamma ,\Bbb D)^0\to
H^\ast (\Gamma_0 ,L_\lambda)^0$$
is $I_{\lambda}H^\ast (\Gamma,\Bbb D)^0$ where $I_{\lambda}$ is the
kernel of the
homomorphism
$\lambda^\dag : \Lambda\lra \z_p$ induced by $\lambda$.
\endproclaim

For easy reference, we will also cite here the
two main lemmas from [Ash-Stevens 97] that we need to
prove Theorem 2.1. We first need a few definitions. We let
$$
W := \spf(\Lambda)
$$
be the formal scheme associated to $\Lambda$.  We note that for each
local field
$K/\Q_p$ the $K$-points of $W$ are in one-one correspondence with
continuous
homomorphisms $k : T(\z_p) \lra K^\times$. If $t\in T(\z_p)$ and $k\in
W(K)$, then
we will write $t^k \in K^\times$ for the value of $k$ on $t$.
Dominant integral weights can (and will) be viewed as elements of
$W(\Q_p)$. We let $W^+$ denote the set of dominant integral weights.

Now let $k\in W(\Q_p)$. We let $C_k(X)$
denote the $\z_p$-module
$$
C_k(X) := \left\{\,f:\X\lra \z_p\,\biggm|\,\matrix \hbox{\rm $f$ is
continuous, and
}\hfill
\\
\hbox{\rm $f(tx) = t^kf(x)$ for
all
$t\in T(\z_p),\ x\in X$}
\endmatrix
\,\right\}.
$$
Let $\bD_k$ be the space of continuous $\z_p$-valued
linear functionals on $C_k(X)$. Clearly
$C_k(X)
\subseteq C(X)$ is a closed subspace, so restriction of linear
functionals
induces a surjective map
$
\bD \lra \bD_k.
$
We regard both $\bD$ and $\bD_k$ as $\Lambda[S]$-modules
(with the $\star$-action of $S$).
Let $I_k\subseteq \Lambda$ be the prime associated to $k\in W(\Q_p)$.
The action of
$\Lambda$ on $\bD_k$ factors through $\Lambda/I_k$.
\medskip

\proclaim{Lemma 2.2}
The ideal $I_k$ is generated by a $\bD$-regular sequence. Moreover, the
canonical sequence
$$
0\lra I_k\bD \lra \bD \lra \bD_k\lra 0
$$
is exact.
\endproclaim

\medskip

\noindent
{\bf Proof:} Lemma 1.1 of [Ash-Stevens] shows that $I_k$ is generated
by a finite
$\bD$-regular sequence. We've already noted above that the
map $\bD\lra\bD_k$ is surjective. We clearly have $I_k\bD$ is contained
in the kernel
and it follows from Lemma 6.3 of [Ash-Stevens 97] that $I_k\bD$
contains the kernel.
Thus the sequence is exact. This completes the proof.
\medskip

\proclaim{Lemma 2.3}
If $\lambda\in W^+$ is a dominant integral weight, then the canonical
map
$$
H^\ast(\Gamma,\bD_\lambda)^0 \buildrel\sim\over \lra
H^\ast(\Gamma_0,L_\lambda)^0
$$
is an $\cH$-isomorphism.
\endproclaim
\medskip

\noindent
{\bf Proof.} The surjectivity of the map follows from Lemma 6.7 and the
injectivity
follows from Lemma 6.8 of [Ash-Stevens 97]. (Note, however, that the
statement of
Lemma 6.8 contains a typo. The conclusion of the statement should read:
``Then $\zeta$ can be
represented by a cocycle $z_m$ in $\Hom_\Gamma(F_k,\bD)$ which takes
values in
$\bD^\ast$.")
\medskip

\noindent
{\bf Remark 2.4.} In an upcoming paper by the first and third authors,
we
generalize the ideas of this section to the non-ordinary ``finite
slope" situation.

\medpagebreak

\heading{Section 3:~~Eigenpackets: Twists and Deformations}
\endheading

\noindent
{\bf Definition 3.0} If $A$  is a $\Lambda$-module,
a $\Lambda$-homomorphism $\varphi :
\cH\lra A$ will be called an $A$-valued eigenpacket.
If $\varphi^\prime:\cH\lra A^\prime$ is another
eigenpacket then a morphism $f: \varphi \lra \varphi^\prime$ is
defined to
be a $\Lambda$-morphism $f : A\lra A^\prime$ such
that the diagram
$$
\matrix
\cH &\buildrel\varphi\over\lra & A\\
\\
&{\varphi^\prime}\atop\searrow&\ \ \big\downarrow f\\
\\
&&A^\prime
\endmatrix
$$
is commutative.
Conversely, if $f : A\lra A^\prime$ is a $\Lambda$-morphism,
then we may push
$\varphi$ forward by $f$ to obtain an eigenpacket $\varphi^\prime :=
f_\ast\varphi$ and
a morphism
$f : \varphi \lra \varphi^\prime$.
If $\varphi$ is an eigenpacket, we will denote its target by 
$A_\varphi$.
\medskip

If $A$ is either a $\Q_p$-Banach algebra or a $\z_p$-algebra that is
separated and complete in the $p$-adic topology, then there is a one-one
correspondence between continuous homomorphisms
$\Lambda \lra  A$ and continuous
characters $T(\z_p) \lra A^\times$.    This is also true for any
$\Q_p$-Banach algebra $A$.   For a character
$\chi : T(\z_p)\lra  A^\times$
we let $\chi^\dagger : \Lambda \lra A$ denote the associated
homomorphism.
Conversely, if $A$ is a $\Lambda$-algebra, we let $\chi_A$ denote the
associated
character and call $\chi_A$ the structure character of $A$. For any
eigenpacket $\varphi$ we define
the {\it weight} of $\varphi$ to be the structure character
$\chi_{A_\varphi}$, which we will denote simply $\chi_\varphi$.

Let $K/\Q_p$ be a finite field extension.  Whenever we say
``$K$-valued eigenpacket'' we will implicitly assume that $K$ has been
endowed with a $\Lambda$-algebra structure.
\medskip

We are interested in eigenpackets that occur in the
$\Gamma_0$-cohomology of finite dimensional irreducible rational
representations of
$GL(n)$. The weight of any such eigenpacket must be a dominant integral
character of
$T(\z_p)$.
\medskip

\noindent
{\bf Definition 3.1} Let $K/\Q_p$ be a finite field extension and let
$\varphi$ be a
$K$-valued eigenpacket. We say $\varphi$ is {\it arithmetic of weight
$\lambda$}, where
$\lambda$ is a dominant integral character,
if
$\varphi$ occurs in
$H^\ast(\Gamma_0,V_\lambda(K))$ (i.e.\ is the homomorphism of
eigenvalues of some Hecke
eigenvector).
\medskip

\medskip

\noindent
{\bf Definition 3.2.}
Given an $A$-valued eigenpacket $\varphi$ and a continuous character
$\psi : \z_p^\times \lra A^\times$ we may ``twist" $\varphi$ by $\psi$
as
follows. Let $A(\psi)$ be the unique continuous
$\Lambda$-algebra whose
underlying
$\z_p$-algebra is $A$ and whose structure character is $\chi_A\cdot
(\psi\circ\det)$.
The $\psi$-twist of $\varphi$ is defined to be
the unique eigenpacket
$$
\tw_\psi(\varphi) : \cH \lra A(\psi)
$$
given on Hecke operators by $U\longmapsto \varphi(U)$ and
$T(\sigma) \longmapsto \psi(\det(\sigma))\varphi(T(\sigma))$ for
$\sigma\in S^\prime$.
In particular, we define the {\it standard twist} of $\varphi$ to be
$$
\tw(\varphi) := \tw_{\bold 1}(\varphi)
$$
where ${\bold 1} : \z_p^\times \lra A^\times$ is the identity character.
Iteration of the standard twist
gives us the {\it $n$-fold twist} $\tw^n(\varphi)$ for any integer $n$.
We also define the {\it universal twist} of $\varphi$ by
$$
{\buildrel\sim\over \tw}(\varphi) := \tw_{\tilde\psi}(i_\ast\varphi)
$$
where $i : A \lra A[[\z_p^\times]] :=
A\otimes_{\z_p}\z_p[[\z_p^\times]]$ is the base
extension morphism and $\tilde\psi : \z_p^\times \lra A[[\z_p^\times]]$
the universal character.
\bigskip

\noindent
Note that by definition, the weight of $\tw_\psi(\varphi)$ is the
product of $\chi_\varphi$ and $\psi\circ\det$.
Moreover, for any $\psi : \z_p^\times \lra A^\times$
there is a unique $A$-morphism ${\buildrel\sim\over\tw}(\varphi) \lra
\tw_\psi(\varphi)$.
The next proposition states that arithmeticity of eigenpackets is
invariant under twist.
\medskip

\noindent
\proclaim{Proposition 3.3}
If $\varphi$ is arithmetic, then $\tw^n(\varphi)$
is arithmetic for every $n\in \z$. If, moreover, $\varphi$ has weight
$\lambda$, then $\tw^n(\varphi)$ has weight $\lambda\cdot\det^n$.
\endproclaim

\noindent
{\bf Proof:} Let $\lambda$ be a dominant integral weight and let
$\varphi$ be an
eigenpacket occurring in $H^\ast(\Gamma_0,V_\lambda(K))$. Now let
$V_\lambda(n)$ be $V_\lambda$ with the $G$-action twisted by the
$n$-fold power
of the determinant. Then the identity map $V_\lambda\lra V_\lambda(n)$
is a $\Gamma_0$-isomorphism and therefore induces an isomorphism
$H^\ast(\Gamma_0,V_\lambda(K)) \lra H^\ast(\Gamma_0,V_\lambda(n)(K))$.
A simple calculation shows that this map sends any $\varphi$-eigenvector
to a $\tw^n(\varphi)$-eigenvector. Thus $\tw^n(\varphi)$ occurs in
$H^\ast(\Gamma_0,V_\lambda(n)(K))$. It follows that $\tw^n(\varphi)$
is arithmetic of weight $\lambda\cdot \det^n$. This completes the proof
of Proposition
3.3.
\bigskip

Consider a topological $\Lambda$-algebra $A$
which is of one of the following two types:
\itemitem{(1)} a $\Q_p$-affinoid algebra; or
\itemitem{(2)} a complete semi-local noetherian $\Lambda$-algebra with
finite residual
fields.
\smallskip

\noindent
In either case, we assume the structure morphism $\Lambda\lra A$ is
continuous.
\medskip

\noindent
{\bf Remark:} Alternatively, we could work in the category of
affine formal schemes over $\spf(\Lambda)$, and indeed for all of
our applications it suffices to work in that category.   In that case,
the ring of rigid analytic functions on an admissible affinoid open set
is of type (1), while on an affine formal scheme finite over
$\spf(\Lambda)$, the ring of global rigid analytic functions of norm
$\le 1$ is of type (2).

If $A$ is arbitrary of type (1) or (2) we will say
$A$ is {\it
$p$-analytic} and  in this case
an $A$-valued eigenpacket will be called $p$-analytic. In case (1) we
will
say $\varphi$ is an {\it affinoid} eigenpacket. These conventions are
justified by the
fact that when $A$ is $p$-analytic there is a naturally defined
$p$-adic rigid analytic space $\cX_A$ over $\Q_p$ such that for every
finite
extension $K/\Q_p$ the $K$-valued points of $\cX_A$ are given by
continuous
homomorphisms from $A$ to $K$:
$$
\cX_A(K) = \Hom_{cont}(A,K).
$$
(For the construction of $\cX_A$ in case (2) see section 7 of [de Jong
96],
or section 1.1 of [Coleman-Mazur 98].)
We will also write $\cX_A(\oQ_p) :=
\Hom_{cont}(A,\oQ_p)$ for the set of
$\oQ_p$-valued points and note that every $\oQ_p$-valued point
must actually take values in some finite extension of $\Q_p$.

We give $\cX_A(\oQ_p)$ the Zariski topology whose closed sets are the zero
sets of the ideals of $A$.  Note that if $A$ is of type (1) this is the
usual topology, but if $A$ is of type (2) it is courser than the usual
topology defined in terms of the full ring of rigid analytic
functions. 
\medskip

Henceforth, all eigenpackets will be assumed to be $p$-analytic. Thus
for every
eigenpacket $\varphi$, the $\Lambda$-algebra $A_\varphi$ is
$p$-analytic.
We denote the associated rigid analytic space by $\cX_\varphi$. We will
often use geometric language to describe properties of $\varphi$. For
example, we say $\varphi$ is affinoid if $A_\varphi$ is an affinoid
algebra and we
say $\varphi$ is irreducible if $\cX_\varphi$ is an irreducible
topological space.
\medskip

\noindent
{\bf Definition 3.4.} A subset $\cS\subseteq \cX_A(\oQ_p)$ is said to
be Zariski dense if $0$ is the only element of $A$ that vanishes at
every
point in $\cS$.
\bigskip

We define the {\it weight space} $\cW$ for
$GL(n)$ by
$\cW := \cX_{\Lambda}$.
For any $p$-analytic $\Lambda$-algebra $A$ the homomorphism
$\Lambda\lra A$
induces a rigid analytic morphism
$$
\wt: \cX_A\lra \cW
$$
which we call the {\it weight map}. We let $\cW^+$
be the set of all dominant integral weights and if $A$ is a $p$-analytic
$\Lambda$-algebra then for any field
$K$ finite over
$\Q_p$ we define
$$
\cX_A^+(K) := \left\{\,x \in \cX_A(K)\,\biggm|\,\wt(x)\in
\cW^+\,\right\}.
$$
We let $\cX_A^+ := \cX_A^+(\oQ_p)$ and call the points of $\cX_A^+$
the dominant integral points on $\cX_A$.
\medskip

\noindent
{\bf Remark 3.5.} The set $\cW^+$ of dominant integral points is a
Zariski dense subset
of $\cW$. However the set $\cX_A^+$ is not, in general, Zariski dense
in $\cX_A$.
\bigskip

If $\varphi$
is an eigenpacket then we say $\varphi$ is an eigenpacket over
$\cX_\varphi$. For
each rigid analytic subspace $Z\subseteq
\cX_\varphi$, we let $\varphi_Z$ be the composition of $\varphi : \cH
\lra A_\varphi$
with the map
$A_\varphi \lra A(Z)$ defined by restriction to $Z$ of rigid analytic
functions on
$\cX_\varphi$. In this case $\varphi_Z$ will be called the restriction
of $\varphi$ to
$Z$. In particular, if $x\in \cX_\varphi(\oQ_p)$, then the restriction
of $\varphi$ to
$\{x\}$ will be called the specialization of $\varphi$ at $x$ and
denoted
$\varphi_x$. We define
$$
\cA(\varphi) := \left\{\,x\in \cX_\varphi(\oQ_p)\,\biggm|\,\varphi_x\
\hbox{\rm is
arithmetic}\,\right\}
$$
and note that $\cA(\varphi)\subseteq \cX_\varphi^+$.
\bigskip

\noindent
{\bf Definition 3.6.}
An eigenpacket $\varphi$ is said to be
\itemitem{(1)} {\it arithmetic} if the set $\cA(\varphi)$ is Zariski
dense
in $\cX_\varphi$, and
\itemitem{(2)} {\it $\oQ_p$-valued} if the underlying ring of $A_\varphi$ is
a finite field extension of $\Q_p$.
\bigskip

If $\varphi$ is an arithmetic eigenpacket, then we regard the family of
all
specializations
$\varphi_x$, $x\in\cX_\varphi(\oQ_p)$, as an analytic family of
$\oQ_p$-valued
eigenpackets that interpolates the arithmetic eigenpackets
$\varphi_x$, $x\in \cA(\varphi)$. We sometimes prefer to take a more
{\it local}
point of view, starting with a $\oQ_p$-eigenpacket and then
``deforming" it to an
analytic family. With this picture in mind, we make the following
definition.
\medskip

\noindent
{\bf Definition 3.7.} Given a $K$-valued eigenpacket $\varphi_0$, a {\it
deformation} of $\varphi_0$ is a pair $(\varphi,x)$ consisting of an
eigenpacket
$\varphi$ and a point $x\in \cX_\varphi(K)$ such that
\itemitem{(1)} $A_\varphi$ is a $p$-analytic integral domain; and
\itemitem{(2)} $\varphi_x = \varphi_0$.

\noindent
We say $(\varphi,x)$ is an arithmetic deformation of $\varphi_0$ if, in
addition to (1)
and (2), we also have
\itemitem{(3)} $\varphi_Z$ is arithmetic for every admissible affinoid
neighborhood
$Z$ of $x$ in $\cX_\varphi$.
\bigskip

Although we won't have use for it later, we record here a notion of equivalence of deformations.  We say two deformations $(\varphi,x),\ (\varphi^\prime,x^\prime)$ of
$\varphi_0$ are locally
equivalent if there exists an irreducible affinoid neighborhood $Z$ of
$x$ in
$\cX_\varphi$, an irreducible affinoid neighborhood $Z^\prime$ of
$x^\prime$ in
$\cX_{\varphi^\prime}$, and an isomorphism $f:
Z^\prime\buildrel\sim\over\lra Z$
such that $f(x^\prime) = x$ and $f^\ast\varphi_Z =
\varphi^\prime_{Z^\prime}$.
\medskip

Let $\tilde\varphi_0 = {\buildrel\sim\over\tw}(\varphi_0) : \cH\lra
\oQ_p[[\z_p^\times]]$
be the universal twist of $\varphi_0$ and let $x_0\in
\cX_{\tilde\varphi_0}(\oQ_p)$
be the point associated to the trivial character of
$\z_p^\times$. It follows from proposition 3.3 that if $\varphi_0$ is
arithmetic, then
the restriction of $\tilde\varphi_0$ to the connected component of
$x_0$ is an arithmetic
deformation of $\varphi_0$ at $x_0$.
\medskip

\noindent
{\bf Definition 3.8.}
An arbitrary deformation $(\varphi,x)$ of $\varphi_0$ will be
called a {\it twist deformation} of $\varphi_0$ if there is a morphism
$f:
\tilde\varphi_0 \lra
\varphi$ such that $f^\#(x) = x_0$.
\bigskip

\noindent
Note that if $(\varphi,x)$ is a twist deformation
of $\varphi_0$, then every specialization $\varphi_z$, $z\in
\cX_\varphi(\oQ_p)$,
is a twist of $\varphi_0$ by some character $\z_p^\times\lra
\oQ_p^\times$.

Clearly, every twist deformation of an arithmetic $\oQ_p$-eigenpacket
$\varphi_0$ is
arithmetic. If the converse is true, we say $\varphi_0$ is
arithmetically rigid.
\medskip

\noindent
{\bf Definition 3.9.} An arithmetic $\oQ_p$-eigenpacket is said to be
{\it arithmetically
rigid} if every arithmetic deformation of $\varphi_0$ is a twist
deformation.
\medskip

In section 9 of this paper, we will give examples of arithmetic
$K$-eigenpackets that are arithmetically rigid.

\medpagebreak

\heading{Section 4:~~The Universal $p$-Ordinary Arithmetic Eigenpacket}
\endheading

In this section we construct the universal $p$-ordinary arithmetic
eigenpacket
$\vpoa$ and use it to give a criterion for rigidity of $p$-ordinary
arithmetic
$\oQ_p$-eigenpackets. Here $K$, as always, denotes a local field finite
over $\Q_p$.
\medskip

\noindent
{\bf Definition 4.1.} A $\oQ_p$-valued eigenpacket $\varphi$ is said to
be {\it
$p$-ordinary} if $\varphi(U)\in \OoQ^\times$. More generally, a
$p$-analytic
eigenpacket
$\varphi$ is said to be $p$-ordinary if $\varphi_x$ is $p$-ordinary for
every $x\in
\cX_\varphi(\oQ_p)$.
\bigskip

We have the following important result.
\bigskip

\noindent
\proclaim{Theorem 4.2} There is a $p$-ordinary arithmetic eigenpacket
$\varphi_{ord}$ satisfying the following universal property: for every
$p$-ordinary arithmetic eigenpacket $\varphi$ there is a unique morphism
$\varphi_{ord}\lra \varphi$.  Moreover, a point $x\in
\cX_{\varphi_{ord}}$ is arithmetic if and only if $x$ has dominant
integral weight.  In other words:
$
\cA(\varphi_{ord}) = \cX_{\varphi_{ord}}^+.
$
\endproclaim

\medskip
Our construction of $\vpoa$ will use the ordinary part of the
cohomology of the universal
highest weight module
$\bD$ constructed in $\S$2. Consider the $\cH$-module $\bH$ defined by
$$
\bH := \bigoplus_i H^i(\Gamma,\bD)^0.
$$
We
define
$\Ro$ to be the image of $\cH$ in $\end_{\Lambda}(\bH)$ and consider the
associated
$\Lambda$-morphism
$$
\vpo : \cH\lra \Ro.
$$
We know from Theorem 2.1 that $\bH$ is a finitely generated
$\Lambda$-module,
and from this it follows that $\Ro$ is a finite $\Lambda$-algebra. Thus
$\vpo$ is a $p$-analytic eigenpacket. Let $\cX_{\vpo}$ be the associated
rigid analytic space, and let $\cA(\vpo)$ be the set of arithmetic
points on $\cX_{\vpo}$ as defined in the last section.
\medskip

\proclaim{Lemma 4.3}
Let $\varphi$ be a $p$-ordinary $\oQ_p$-valued eigenpacket. Then
$\varphi$ is arithmetic if and only if $\varphi = (\vpo)_x$ for
some $x\in\cX_{\vpo}^+$.  In particular,
$\cA(\vpo) =
\cX_{\vpo}^+.
$

\endproclaim

The proof of Lemma 4.3 will be given in the next section. In this
section, we
assume Lemma 4.3 and prove Theorem 4.2 as a consequence.
\medskip

\noindent
{\bf Proof of Theorem 4.2.} Let $I_{arith}\subseteq \Ro$
be the ideal defined by
$$
I_{arith} := \bigcap_{x\in\cA(\vpo)}\ker(x).
$$
  We then let $\Roa := \Ro/I_{arith}$ and define
$$
\vpoa : \cH \lra \Roa
$$
to be the composition of $\vpo$ with $\Ro\lra \Roa$. By definition, we
see that
$\vpoa$ is arithmetic. Let $J := \ker(\vpoa) \subseteq \cH$.

Now let $\varphi : \cH \lra A_\varphi$ be an arbitrary $p$-ordinary
arithmetic
eigenpacket. We will show $J\subseteq \ker(\varphi_x)$ for every
$x\in\cA(\varphi)$.
Indeed, if $x\in \cA(\varphi)$ then, since $\varphi_x$ is arithmetic,
Lemma 4.3 tells us that there is a point $y\in \cA(\vpo) = \cX_{\vpo}^+$
such that
$\varphi_x = (\vpo)_y$. But we clearly have $J\subseteq
\ker((\vpo)_y)$, so
$J\subseteq \ker(\varphi_x)$, as claimed.

It now follows easily that $J\subseteq \ker(\varphi)$. Indeed, from the
last
paragraph we have
$$
\varphi(J)\subseteq \bigcap_{x\in \cA(\varphi)}\ker(x),
$$
and the intersection on the right is $\{0\}$ since $\varphi$ is
arithmetic.
Thus we have shown that $\ker(\vpoa)\subseteq \ker(\varphi)$. Since
$\vpoa : \cH \lra \Roa$ is surjective by definition, we conclude that
$\varphi : \cH \lra A_\varphi$ factors through $\vpoa$. This proves the
existence of a morphism $\vpoa \lra \varphi$. Uniqueness follows from
the surjectivity
of $\vpoa$. This completes the proof of Theorem 4.2.
\bigskip

It is natural to ask whether or not the morphism $\vpo\lra\vpoa$ is an
isomorphism,
or equivalently, whether $\vpo$ is arithmetic.
Hida's theory answers
this in the affirmative when $n=2$, but we know almost nothing in this
direction when
$n\ge 3$. The eigenpacket $\vpo$ deserves to be of interest independent
of whether
or not it is arithmetic. For example, $\vpo$ parametrizes $p$-ordinary
torsion
in the cohomology of the $\z_p$-lattices $L_\lambda$, $\lambda\in
\cW^+$.
If $\vpo \lra\vpoa$ is an isomorphism, then it can be shown that every
eigenpacket
occurring in the $p$-ordinary torsion must be congruent to some
$p$-ordinary arithmetic
$\oQ_p$-eigenpacket.

\noindent
\proclaim{Definition 4.4} We define the {\it $p$-ordinary eigenscheme}
to be the
formal scheme $\Xoa := \spf(\Roa)$. We also let $\cXoa := \cX_{\vpoa}$
and call this the $p$-ordinary eigenvariety.
\endproclaim

If $\OK$ is the ring of integers in a finite extension $K/\Q_p$ then
any continuous
homomorphism $\Roa\lra K$ takes values in $\OK$. Hence the canonical map
$$
X_{ord}(\OK) \buildrel\sim\over\lra \cXoa(K)
$$
is a bijection. Thus the scheme $X_{ord}$ remembers the $K$-points of
$\cXoa$.
But the scheme $X_{ord}$ contains additional information related to
congruences
and this will be important to us below in the definition of quasi-cuspidality. 
For $x\in X_{ord}(K)$ we let $P_x$
be the kernel in
$\Roa$ of
$x$. Then $P_x$ is a prime ideal in $\Roa$.
\medskip

\noindent
{\bf Definition 4.5.} If $x,x^\prime\in \cXoa(\oQ_p)$ are two
$\oQ_p$-points, then we
say
$x$ and
$x^\prime$ are {\it congruent} if $P_x$ and $P_{x^\prime}$ are contained
in the same maximal ideal of $\Roa$.
\medskip

We note that two points $x,x^\prime\in \cXoa(\oQ_p)$ are congruent if
and only if
there is an automorphism $\tau$ of $\oQ_p/\Q_p$ such that for every
$\alpha\in \cH$
we have the congruence
$(\vpoa)_x(\alpha) \equiv \tau\bigl((\vpoa)_{x^\prime}(\alpha)\bigr)$
modulo the
maximal ideal of the integers of $\oQ_p$. This is the usual notion of
congruence
between two systems of Hecke eigenvalues.
Indeed, since $\Roa$ is
complete and semilocal,
the set of maximal ideals of
$\Roa$ is in one-one correspondence with the set of connected
components of
$\Xoa$ and these are in one-one correspondence with the set of
congruence classes of
points on $\cXoa$.

\medskip

We close this section by showing that $\vpoa$ is isomorphic to
$\tw(\vpoa)$ and
deriving some simple consequences.
\medskip

\noindent
{\bf Theorem 4.6.} {\it
The unique morphism $\vpoa \lra
\tw(\vpoa)$ is an isomorphism.  Moreover,
for every specialization $\varphi_0$ of $\vpoa$ at a point of
$\cXoa(\oQ_p)$,
there is a morphism
$$
\vpoa \lra \tws(\varphi_0)
$$
where $\tws(\varphi_0)$ is the universal twist of $\varphi_0$ defined
in section 3.
}
\medskip

\noindent
{\bf Proof:} Recall that if $A$ is a $\Lambda$-module, then $A({\bold 
1})$
denotes its twist by $\det$ (Definition 3.2).
Let $\bD \cong \bD(X)$ be the $\Lambda[S]$-module
described in section 2.
We define the $\Lambda$-morphism $\tw : \bD \lra \bD({\bold 1})$
by sending a measure $\mu\in \bD$ to the measure $\tw(\mu)\in
\bD({\bold 1})$
defined by the integration formula
$$
\int_X f(x)d\tw(\mu)(x) = \int_X \det(x)f(x)d\mu(x)
$$
for $f\in C(X)$. Clearly $\tw$ is an isomorphism. On the other hand,
for $\sigma\in S$
we have $\tw(\mu|\sigma) = \det(\sigma)\tw(\mu)|\sigma$. In particular,
$\tw$ is a
$\Gamma$-morphism and therefore induces a $\Lambda$-isomorphism
$$
\bH \buildrel\sim\over\lra \bH({\bold 1}).
$$
Conjugation induces a $\Lambda$-isomorphism
$
\tw : \end_\Lambda(\bH) \buildrel\sim\over\lra \end_\Lambda(\bH({\bold 
1})).
$
A computation shows that the diagram
$$
\matrix
\cH &\buildrel\tw\over\lra & \cH({\bold 1}) \\
\\
\big\downarrow && \big\downarrow \\
\\
\end_\Lambda(\bH) &\buildrel\tw\over\lra &\end_\Lambda(\bH({\bold 1}))
\endmatrix
$$
is commutative where the vertical arrows are both given by $\vpo$.
Hence $\tw$ induces an isomorphism $\tw: \vpo \lra \vpo$. By
Proposition 3.3,
this morphism preserves the ideal $I_{arith}$ and therefore induces an
isomorphism
$\tw :\vpoa \lra \tw(\vpoa)$.

Now let $\varphi_0$ be the specialization of $\vpoa$ at a point $x_0\in
\cXoa(\oQ_p)$.
For each integer $n$, let $x_n:=\tw^n(x_0)\in \cXoa(\oQ_p)$ and note
that we
then have a commutative diagram
$$\matrix
\cH &\buildrel\vpoa\over\lra & \Roa \\
\big\downarrow && \big\downarrow \\
\oQ_p[[\z_p^\times]] & \buildrel\psi_n^\dagger\over\lra & \oQ_p \\
\endmatrix
$$
where the left vertical arrow is $\tws(\varphi_0)$, the right vertical
arrow is $x_n$, and $\psi_n^\dagger$ is the map whose structure
character
is the $n$th power character $\z_p^\times \lra \oQ_p^\times$. But since
$\displaystyle\bigcap_n \ker(\psi_n^\dagger) = \{0\}$, we see that
$\tws(\varphi_0)$ maps every $\alpha\in \ker(\vpoa)$ to $0$. Finally,
since
$\vpoa$ is surjective, we conclude that
$\tws(\varphi_0)$ factors through $\vpoa$. This gives us a morphism
$\vpoa \lra \tws(\varphi_0)$. This completes the proof of Theorem 4.6.

\bigskip

If $\varphi=(\vpoa)_x$ is the eigenpacket associated to an arbitrary 
point
$x\in\cXoa(\oQ_p)$, then the last theorem gives us a morphism
$\Roa \lra \overline \Q_p[[\z_p^\times]]$. If we compose this morphism
with the morphism associated to the canonical projection $\z_p^\times
\lra 1+p\z_p$ (or to $1+4\z_p$ if $p=2$)
we obtain a $\Lambda$-morphism
$$
\tilde x : \Roa \lra \oQ_p[[1+p\z_p]]
$$
whose kernel is a prime ideal $P_{\tilde x}\subseteq \Roa$.
\medskip

\noindent
{\bf Definition 4.7.} For arbitrary $x\in \cXoa(\oQ_p)$ we let
$Z(\tilde x) \subseteq \cXoa$ be the irreducible Zariski closed subset
of
$\cXoa$ cut out by $P_{\tilde x}$.
\medskip

We have the following corollary of Theorem 4.6.

\proclaim{Corollary 4.8}
If $Z$ is an irreducible component of $\cXoa$ and $x\in \cXoa(\oQ_p)$
then
$x\in Z(\oQ_p)$ $\Longrightarrow$ $Z(\tilde x)\subseteq Z$.
\endproclaim

\noindent
{\bf Proof:} Since $\Roa$ is finite over $\Lambda$, we have $\cXoa$ is a
noetherian topological space (with the Zariski topology), hence has
only finitely many
irreducible components. Consider the powers $\tw^n$ of the twist
morphism $\tw : \cXoa
\lra
\cXoa$ of Theorem 4.6 for $n\equiv 0$ modulo $(p-1)$ (or modulo $2$ if $p=2$). Clearly,
$\tw^n(Z)$ is
an irreducible component of $\cXoa$ for every such $n$, thus there is
an integer
$N\equiv 0$ modulo $(p-1)$ (again modulo $2$ if $p=2$) such that $\tw^n(Z) = \tw^m(Z)$ whenever $n
\equiv m$ modulo
$N$. In particular, we have $\tw^n(Z) = Z$ whenever $N|n$. From this it
follows
that the set $\{\,\tw^n(x)\,|\,N|n\,\} \subseteq Z(\oQ_p)$. But this
set is a Zariski
dense subset of $Z(\tilde x)$. Therefore $Z(\tilde x) \subseteq Z$ and
Corollary
4.8 is proved.

\proclaim{Theorem 4.9}
Let $\varphi_0$ be a $p$-ordinary arithmetic $\oQ_p$-eigenpacket
and let $x_0\in \cXoa^+$ be the corresponding point.
Then $\varphi_{0}$ is rigid if and only if $Z(\tilde x_0)$ is
an irreducible component of $\cXoa$.
\endproclaim

\noindent
{\bf Proof.} If $\varphi_0$ is rigid and $Z$ is an irreducible
component of $\cXoa$
containing $x_0$, then by the definition of rigidity, $Z$ must be a
twist-deformation
of $x_0$. Hence every arithmetic point of $Z$ is a twist of $x_0$ and
therefore
every arithmetic point of $Z$ is in $Z(\tilde x_0)$. Since the
arithmetic points of
$Z$ are Zariski dense in $Z$, we conclude that $Z\subseteq Z(\tilde
x_0)$.
So by Corollary 4.8, we have $Z = Z(\tilde x_0)$ as claimed.

Conversely, suppose $Z(\tilde x_0)$ is an irreducible component of
$\cXoa$
and let $(\varphi,z_0)$ be any arithmetic deformation of $\varphi_0$.
Since $\varphi_0$
is ordinary, we can choose an irreducible affinoid neighborhood $Z$ of
$z_0$ in
$\cX_\varphi$ such that $\varphi_Z$ is $p$-ordinary. Hence there is a
morphism
$f: \vpoa \lra \varphi_Z$ such that $f^\#(z_0) = x_0$. Since $Z$ is
irreducible,
we have also $f^\#(Z)$ is irreducible in $\cXoa$. Let $W$ be an
irreducible component containing $f^\#(Z)$. By Corollary 4.8, $Z(\tilde
x_0)\subset W.$   Since $Z(\tilde x_0)$ is an irreducible component,
$Z(\tilde x_0) = W$. Thus $f^\#(Z)
\subseteq Z(\tilde
x_0)$. Hence $\varphi_Z$ is a twist of $\varphi_0$ and it follows that
$(\varphi,z_0)$
is a twist deformation of $\varphi_0$. This proves $\varphi_0$ is rigid
and Theorem 4.9
is proved.

\medpagebreak

\heading{Section 5:~~Proof of Lemma 4.3.}
\endheading

In this section, we prove a general cohomological result (Theorem 5.1)
and deduce
Lemma 4.3 as a consequence.

Let
$M$ be a profinite
$\Lambda[S]$-module and let
$H^\ast(M)^0$ be the ordinary part of the cohomology $H^\ast(\Gamma,M)$
with
respect to the operator
$T(\pi)$. Then we may (and do) regard
$H^\ast(M)^0$ as an
$\cH$-module. Define
$R(M)$ to be the image of
$\cH$ in
$\end_{\Lambda}(H^\ast(M)^0)$ and define $A(M)$ to be $R(M)^{red}$, i.e.
$R(M)$ modulo its nilradical.
\bigskip

\proclaim{Theorem 5.1} Let $M$ be a profinite $\Lambda[S_0]$-module and
suppose
$H^\ast(M)^0$ is a finitely generated $\Lambda$-module.
Let $I$ be an ideal in $\Lambda$ that is generated by a finite
$M$-regular sequence. Then $H^\ast(M/IM)^0$ is a finitely generated
$\Lambda/I$-module. Moreover there is a
natural isomorphism
$$
R(M)/Rad_{R(M)}(I) \cong A(M/I)
$$
where $Rad_{R(M)}(I)$ is the radical of $IR(M)$ in $R(M)$.
\endproclaim

Before proving Theorem 5.1, we first indicate how it implies Lemma 4.3.
\medskip

\noindent
{\bf Proof of Lemma 4.3.} If we take $M =\bD$ then we have
$R(\bD) = \Ro$. Moreover, we know (see Lemma 2.2) that for each
dominant integral weight $\lambda$, the ideal
$I_\lambda\subseteq
\Lambda$ is generated by a
$\bD$-regular sequence and that $\bD/I_\lambda \bD \cong \bD_\lambda$.
Hence, by Theorem 5.1, we have
$
A(\bD_\lambda) \cong \Ro/Rad_{\Ro}(I_\lambda).
$
Moreover, from Lemma 2.3 we have
$
A(L_\lambda) \cong A(\bD_\lambda).
$
Thus we have
$$
A(L_\lambda) \cong \Ro/Rad_\Ro(I_\lambda)
$$
from which Lemma 4.3 is an immediate consequence.
\medskip

We now turn to the proof of Theorem 5.1.
\medskip

\noindent
{\bf Proof of Theorem 5.1.} Our proof is by induction on the length of
an
$M$-regular sequence of generators for $I$. We begin with the case
where $I$ is principal, generated by an $M$-regular element $\alpha\in
\Lambda$, i.e.\
an element
$\alpha$ such that $M$ has no $\alpha$-torsion. In this case, we have
an exact sequence
$$
0\lra M \buildrel\alpha\over\lra M\lra M/\alpha M\lra 0.
$$
Now pass to the long exact cohomology sequence, we obtain an exact
sequence
$$
0\lra H/\alpha H \lra H^\ast(M/\alpha M)^0\lra H[\alpha]\lra 0
$$
where $H := H^\ast(M)^0$ and $H[\alpha]$ is the $\alpha$-torsion in
$H$. We first
construct a homomorphism $R(M)\lra A(M/\alpha M)$. For this, we let
$x\in \cH$ be such that $x$ annihilates $H$. Then the last exact
sequence
implies $x^2$ annihilates $H^\ast(M/\alpha M)^0$. Hence $x$ maps to an
element of the
nilradical of $R(M/\alpha M)$ and it follows that $x$ maps to $0$ in
$A(M/\alpha M)$. Thus the canonical map $\cH \lra A(M/\alpha M)$ factors
through the canonical map $\cH\lra R(M)$. Hence we have a
canonical surjective map
$$
\tilde\varphi : R(M)\lra A(M/\alpha M).
$$
We need to show that $\ker(\tilde\varphi) = Rad_{R(M)}(\alpha)$. The
inclusion
$\supseteq$ is obvious.

So let $x\in \ker(\tilde\varphi)$. From the above exact sequence, we
conclude that
$xH\subseteq \alpha H$.
Since, by hypothesis, $H$ is finitely generated over $\Lambda$, there
is a positive
integer
$m$ such that
$H[\alpha^{m+1}] = H[\alpha^m]$. In particular, we
see that $\alpha^m H$ has no $\alpha$-torsion and from this it follows
that
$\alpha : \alpha^mH\lra \alpha^{m+1}H$ is an isomorphism of
$\Lambda$-modules.
Let $\beta : \alpha^{m+1}H\lra \alpha^mH$ be the inverse map. Let
$y\in \end_\Lambda(H)$ be the composition
$$
y: H \buildrel x^{m+1} \over \lra \alpha^{m+1}H \buildrel \beta\over
\lra \alpha^mH
\subseteq H.
$$
Then $\alpha y = x^{m+1}$.

Now define
$$
R^\prime := \{\,y\in\end_\Lambda(H)\,|\, \exists k\in \z^+, z\in R(M)\
\hbox{\rm
such that
$\alpha^ky =\alpha^mz$}\,\}.
$$
The endomorphism $y$ constructed in the last paragraph is an element of
$R^\prime$.
Clearly, $R^\prime$ is a finite $\Lambda$-algebra containing $R(M)$.
Moreover,
for every element of $\rho\in R^\prime$ we have $\alpha^k\rho\in R(M)$
for some
$k\ge 0$. Since $R^\prime$ is finitely generated over $\Lambda$, there
is an
exponent $N$ such that $\alpha^N R^\prime\subseteq R(M)$.

Now consider $x^{(m+1)(N+1)} = (\alpha y)^{N+1} =
\alpha(\alpha^Ny^{N+1})$. Since
$y^{N+1}\in R^\prime$ we have $z:=\alpha^Ny^{N+1}\in R(M)$. Hence
$x^{(m+1)(N+1)}\in
\alpha R(M)$. This proves $x\in Rad_{R(M)}(\alpha)$. Hence
$\ker(\varphi)\subseteq Rad_{R(M)}(\alpha)$. This completes the proof
in the
special case where $I$ is generated by a single $M$-regular element of
$\Lambda$.

Now we suppose $r\ge 1$ and that the theorem is true whenever $I$ is
generated by an
$M$-regular sequence of length $r$. Let
$\alpha_1,\ldots,\alpha_r, \alpha$ be an $M$-regular sequence of length
$r+1$. Let $I$
be the ideal generated by $\alpha_1,\ldots,\alpha_r$ in $\Lambda$ and
let $J$ be the
ideal generated by $I$ and $\alpha$.
We define $\psi$ to be the composition of the surjective homomorphisms
$$
\psi : R(M) \buildrel\varphi\over\lra R(M/IM) \lra A(M/JM).
$$
Let $x\in \ker(\psi)$. Then by the principal case proved in the last
paragraph, we
have $y:=
\varphi(x)\in Rad_{R(M/IM)}(\alpha)$. Thus there is an $m\in \z^+$ such
that $y^m\in
\alpha R(M/IM)$. This means
$\varphi(x^m)\in
\alpha R(M/IM)$. So there is an element $z\in R(M)$ such that
$\varphi(x^m) =
\alpha\varphi(z)$. It then follows from the induction hypothesis that
$x^m-\alpha z \in
Rad_{R(M)}(I)$. There is therefore a positive integer $N$
such that
$$
(x^m -\alpha z)^N\in IR(M).
$$
  From this we see at once that $x^{mN}\in JR(M)$. Hence $x\in
Rad_{R(M)}(J)$.
This proves $\ker(\psi)\subseteq Rad_{R(M)}(J)$. But the opposite
inclusion is
immediate.
Hence $\ker(\psi) = Rad_{R(M)}(J)$. This completes the proof of Theorem
5.1.

\medpagebreak

\heading{Section 6:~~Cuspidal and Quasi-cuspidal Eigenpackets}
\endheading

In this section we recall briefly the connection between cohomology and
automorphic
representations and review some results from the theory of mod $p$
cohomology of
congruence subgroups $\Gamma$ of $SL(n,\z)$ in order to have a
sufficient condition for
arithmetic eigenpackets to be cuspidal automorphic.

We denote the adeles of $\Q$ by $\A$.
Let $\Gamma = \Gamma_0(N)$ for some $N$, $M$ a compact open subgroup of
$GL(n,\A_f)$ whose intersection
with $GL(n,\Q)$ is $\Gamma$. Let $\tilde\frak g$ denote the Lie algebra
of $SL(n,\R)$, $K = O(n)$, $Z$ the center of $SL(n,\R)$. Let $E$ be any
finite
dimensional irreducible complex representation
of $GL(n,\R)$.

The cuspidal cohomology $H_{cusp}^\ast(\Gamma,E)$ is defined in [Borel
83]. By Lemma 3.15, page 121 of [Clozel 90],
there is a canonical surjective map
$$
\Psi: \bigoplus_{\Pi} \
H^\ast(\tilde\frak g, K; \Pi_{\infty} \otimes E) \otimes \Pi_f^M
\to H_{cusp}^\ast(\Gamma,E)
$$
where $\Pi$ runs over a certain finite set of cuspidal automorphic
representations of $GL(n,\A)$.

If $\alpha \in H_{cusp}^\ast(\Gamma,E)$ is an $\cH$-eigenclass, then by
the strong multiplicity one theorem for $GL(n)$,
there is a unique $\Pi_\alpha$ such that
$\alpha \in \Psi(H^\ast(\tilde\frak g, K; \Pi_{\alpha,\infty} \otimes E)
\otimes \Pi_{\alpha,f}^M)$.

We fix an isomorphism of $\C_p$ and $\C$.
Recall that $\Gamma_0$ denotes the subset of $\Gamma$ consisting of
those matrices that are upper triangular modulo $p$.
If $\alpha \in H^\ast (\Gamma_0 ,L_\lambda)$, we say that $\alpha$ is
cuspidal if
its image in $H^\ast (\Gamma_0 ,L_\lambda\otimes\C) \simeq
H^\ast (\Gamma_0 ,L_\lambda) \otimes_{\z_p} \C$ is cuspidal.

\proclaim{Definition 6.0} An $\cO_K$-valued eigenpacket is said to be cuspidal if it occurs as the system of $\cH$-eignevalues on a cuspidal cohomology class.  A point on $\cXoa(K)$ is cuspidal if $\varphi_x$ is cuspidal.
 \endproclaim
 
We will need the following fact: set $n=3$ and
let $\partial Y$ denote
the boundary of the Borel-Serre compactification of $\Gamma\backslash
GL(n,\R)/KZ$. Then
$$
H^3_{cusp}(\Gamma,E) \oplus H^3(\partial Y,E) \simeq H^3(\Gamma,E).
$$
See for instance Proposition 2.1 and its proof in [Ash-Grayson-Green
84] which,
though stated there for trivial coefficients, remains valid for
$E$-coefficients.
\medskip

\medpagebreak

If $\varphi$ is an $\cO_K$-valued eigenpacket such that
$\varphi(T(\ell,k)) = a(\ell,k)$
then the Hecke polynomial of $\varphi$ at $\ell$ is defined to be
$$
P_{\varphi,\ell}(X) = \sum_{k=0}^n (-1)^k a(\ell,k) \ell^{k(k-1)/2}
X^k\in \cO_K[X].
$$
If $x\in \cXoa(K)$ then we will write $P_{x,\ell}$ for
$P_{\varphi_x,\ell}$.
\medskip

\proclaim{Definition 6.1} An $\cO_K$-valued eigenpacket $\varphi$ is
said to be
$\ell$-quasicuspidal if $P_{\varphi,\ell}$ is irreducible modulo the
maximal ideal
of $\cO_K$. We say $\varphi$ is quasicuspidal if $\varphi$ is
$\ell$-quasicuspidal for
some prime $\ell\ne p$.

We say a point $x\in \cXoa(K)$ is
$\ell$-quasicuspidal (resp.
quasicuspidal) if $\varphi_x$ is $\ell$-quasicuspidal (resp.\
quasicuspidal).
\endproclaim
\medskip

The following proposition follows immediately from the definitions.
\medskip

\proclaim{Proposition 6.2} If $x\in \cXoa(K)$ is quasicuspidal
then every point in the connected component of $x$ in $\cXoa(K)$ is also
quasicuspidal.
\endproclaim

\proclaim{Conjecture 6.3} Every quasicuspidal arithmetic eigenpacket is
cuspidal.
\endproclaim
\medskip

\proclaim{Theorem 6.4} If $n=2$ or $3$ then
conjecture 6.3 is true.
\endproclaim

\demo{Proof} This follows from Proposition 3.2.1 of [Ash-Stevens 86]
which says that a Hecke eigenclass in $H^3(\partial Y,E)$ has all its
Hecke polynomials reducible
plus the fact that a non-cuspidal Hecke eigenclass must
restrict nontrivially to the Borel-Serre boundary, if $n \le 3$.
\enddemo

\heading{Section 7: Removing $p$ from the level} \endheading In this
section only, in order to conform with the convention of the
references, all representation spaces will be left modules.

Let $I$ denote the Iwahori subgroup of $GL(n,\z_p)$ consisting of all
matrices which are upper triangular modulo $p$. When $n = 3$, our
computers are not powerful enough to compute the cohomology of
$\Gamma_0 := \Gamma \cap I$. We can only compute for the larger group
$\Gamma$. Therefore we need a theorem that allows us to remove $p$ from
the level of a cohomology class.

Let $G=GL(n)$ and let $(B,T)$ be the pair of upper triangular and
diagonal matrices respectively. We use the corresponding ordering on
the roots, so that the positive roots are those
occurring in $B$. Let $\delta$ be the character on $T$ which is the
$p$-adic
absolute value of the product of the positive root characters:
$\delta(\diag(t_1,\dots,t_n)) =
|t_1|^{n-1}|t_2|^{n-3}\cdots|t_n|^{1-n}$.

Recall that we have identified the Hecke algebras $\cH$ and $\cH_0$ for $\Gamma$ and $\Gamma_0$, and so $\cH$ acts on the cohomology of both $\Gamma$ and $\Gamma_0$.  We then have

\proclaim{Theorem 7.1} Let $\Gamma$ be a congruence subgroup of
$GL(n,\z)$ of
level prime to $p$ and let $\Gamma_0 = \Gamma \cap I$ as above. Let
$\lambda$
be the highest weight $(a_1,a_2,\dots,a_n)$, so that $a_1 \ge a_2 \ge
\dots \ge
a_n$. Let $\phi$ be a packet of Hecke eigenvalues (for $\cH$) that
occurs in
$H^\ast_{cusp} (\Gamma_0 ,L_\lambda(\C))^0$. Suppose $a_1 > a_2 > \dots
  > a_n$.
Then there is an $\cH$-eigenclass $\alpha\in H^\ast_{cusp} (\Gamma ,L_\lambda(\C))$ such that for 
any $l\neq p$ and any $i=1,\ldots,n$, $T(l,i)\alpha=\phi(T(l,i))\alpha.$
\endproclaim

\demo{Proof} Let $|\ |$ denote the $p$-adic absolute value normalized
so that
$|p| = p\inv$.  Let $\alpha \in H^\ast_{cusp} (\Gamma_0 ,L_\lambda)^0$
be an
ordinary cuspidal  cohomology class which is a Hecke eigenclass with
eigenvalues given by $\phi$. Let $\Pi$ be the irreducible automorphic
representation of $GL(n,\A)$ associated to $\alpha$ and write $\Pi_f =
\otimes
\Pi_v$.

Then $\Pi_p$ has a nonzero $I$-invariant vector and therefore is
isomorphic to
a subrepresentation of the induced representation $\I(\chi)$ from
$B(\Q_p)$ to
$G(\Q_p)$ of some unramified character $\chi$ on $T(\Q_p)$ (Theorem 3.8
in
[Cartier 79]). This is normalized induction, so that $\I(\chi)$
consists of the
locally constant functions $f:G(\Q_p) \to \C$ such that $f(bg) =
\delta(b)^{1/2}\chi(b)f(g)$ for $b \in B$ and $g \in G$ (where $\chi$ is
extended from $T$ to $B$ by making it trivial on the unipotent radical
of $B$.)

By a straightforward but rather complicated explicit computation in the 
induced representation, we can determine
the
unramified character $\chi$ given the eigenvalues of the Hecke
operators at
$p$. Let $U(p,m)$ denote the Hecke operator $U(1,\dots,1,p,\dots,p)$
with $m$
copies of $p$, $m = 0, \dots, n$. We may choose $\alpha$ so that it is
also an eigenclass for all the $U(p,m)$ operators. Let $\beta_m$ be the 
eigenvalue of
$U(p,m)$
on the cohomology class $\alpha$ with respect to the $\star$-action. By
assumption, the algebraic integers $\beta_i$ are $p$-adic units. Note
that
$\beta_0 = \beta_n = 1$, since $U(p,0)$ and $U(p,n)$ act via the
$\star$-action
as the identity.

To compare the cohomology with the automorphic theory, we have to
remove the
$\star$-action. With respect to the usual action, the eigenvalue of
$U(p,m)$ on
$\alpha$ is $p^{{a_n}+a_{n-1}+\dots+a_{n-m+1}}\beta_m$.

Let $d_m$ be the diagonal matrix with 1's along the diagonal except for
a $p$
at the $m$-th place. The explicit computation alluded to above gives the
following values, for $m = 1,\dots,n$: $$
\chi\delta^{1/2}(d_m)=p^{a_{n-m+1}-n+2m-1}\beta_m\beta_{m-1}\inv. $$

It follows that if $r_m = {a_{n-m+1}-(n-2m+1)/2}$, then

$$ \chi(d_m)=p^{r_m}\beta_m\beta_{m-1}\inv. $$

Define characters $\chi_i$ of $\Q_p$ by $\chi(\diag(t_1,\dots,t_n)) =
\chi_1(t_1)\dots\chi_n(t_n)$. Since $\chi$ is unramified, each
$\chi_i(t) =
|t|^{c_i}$ for some complex number $c_i$.

It is a result of Bernstein and Zelevinsky ([Bernstein-Zelevinsky 77],
theorem
4.2) that $\I(\chi)$ is irreducible unless $\chi_i\chi_j\inv=|\ |$ for
some
$i\neq j$. See section 1 of [Kudla 94] for a discussion of question of
reducibility of represenations induced from parabolic subgroups.

We compute
$$
\chi_i\chi_j\inv(p) =
p^{r_i-r_j}\beta_i\beta_{i-1}\inv\beta_j\inv\beta_{j-1}.
$$

Note that since $a_1 > a_2 > \dots > a_n$, the $r_m$ are a strictly
increasing
sequence of integers with gaps of at least 2 between successive terms.
Taking
into account the fact that the $\beta_m$ are $p$-adic units, we see that
$\chi_i\chi_j\inv(p)$ never equals $p\inv=|p|$. So we conclude that
$\I(\chi)$
is irreducible and so $\Pi_p=\I(\chi)$.

Let $\Psi$ be the map from $(\tilde\frak g,K)$ cohomology discussed in
section
6. As in that section, $M$ denotes a compact subgroup of $G(\A_f)$ whose
intersection with $G(\Q)$ is $\Gamma$. We let $M_0 = \{y \in M \ | \
y_p \in I
\}$. Then $\alpha = \Psi(c \otimes x)$ for some generator $c$ of
$H^\ast(\tilde\frak g, K; \Pi_{\alpha,\infty} \otimes E)$ and for some
$x =
\otimes x_v \in \Pi_f^{M_0}$.

It follows from (c) p. 152 of [Cartier 79] that $\Pi_p$ contains a
vector
$x_p'$ fixed under $GL(n,\z_p)$. (So in fact $\Pi$ is unramified at
$p$.) Set
$x' = \otimes x'_v$ where $x'_v=x_v$ if $v \ne p$. Then $\alpha'=\Psi(c
\otimes
x')$ is a Hecke eigenclass in $H^\ast_{cusp} (\Gamma ,L_\lambda(\C))$
with the
same packet of $\cH'$-eigenvalues as $\alpha$. \enddemo

\medpagebreak
\heading{Section 8:~~Essential Self-duality}
\endheading

We begin this section with some considerations of selfduality.  If $B$
is a
   ring and $\pi$ a representation of $GL(n,B)$, the contragredient of
$\pi$ is
the representation $\pi'$ defined by $\pi'(g) = \pi({}^tg\inv)$.  We
say that
$\pi$ is ``essentially selfdual" if there is a character $\chi$ of
$B^\times$
such that $\pi'$ is isomorphic to
$\pi \otimes (\chi \circ \det)$.

In parallel to this, if $H$ is a Hecke algebra over $B$ with involution
$\iota$
and $\phi: H \to A$ is a $B$-homorphism, we will say that $\phi$ is
``essentially selfdual" if there is a character $\chi$ of $A^\times$
such that
for every double coset $T$ in $H$,
$\phi(\iota(T)) = \chi(\det(T)) \phi(T)$.

\medskip

\noindent
{\bf Definition 8.1.} a) If $(\Gamma,S)$ is a Hecke pair for $GL(n)$
and $E$
is a finite dimensional $S$-module over a ring of characteristic 0, a
Hecke
eigenclass in $H^*(\Gamma, E)$ is called an essentially selfdual
cohomology
class if its associated automorphic representation (as in $\S 6$) is an
essentially selfdual representation.

b) If $\phi$ is an $\cH$-eigenpacket (as in $\S 3$), it is called an
essentially selfdual eigenpacket if $\phi | \cH'$ is essentially
selfdual,
where $\cH'$ is the Hecke algebra of the semigroup $S'$ defined in $\S
1$ and
$\iota$ is the involution of it given by transpose-inverse composed with
conjugation by the matrix diag$(N,1,\dots,1)$.

\medskip

It follows from the definitions that an essentially
selfdual cohomology eigenclass affords an essentially
selfdual eigenpacket.

\medskip

By theorem 6.12 II in [Borel-Wallach 00], the cuspidal cohomology of an
arithmetic group with coefficients in a finite dimensional irreducible
rational
module $E$ over a field $L$ of characteristic 0 vanishes unless $E$ is
essentially selfdual as a $L$-module. We shall call such $E$ an
essentially
selfdual coefficient module and its highest weight an essentially
selfdual weight.

It's easy to see that an essentially
selfdual eigenpacket has an essentially
selfdual weight.  The converse is not true, as shown by the examples
computed
in [Ash-Grayson-Green 84].
\medskip
The (contravariant) dominant integral weights with respect to
$(B^{opp},T)$ are
given by
$n$-tuples of integers
$\lambda = (\lambda_1,\lambda_2,\ldots,\lambda_n)$ with
$\lambda_1
\ge \lambda_2 \ge \cdots \ge \lambda_n$. Such an $n$-tuple corresponds
to the
character
$t:=\diag(t_1,t_2,\ldots,t_n)
\mapsto t^\lambda :=t_1^{\lambda_1}t_2^{\lambda_2}\cdots
t_n^{\lambda_n}$ on
the torus
$T$. The highest weight of the contragredient of $V_\lambda$ is the
dominant
integral weight $\lambda^\vee := (-\lambda_n,\ldots,-\lambda_1)$. So
$\lambda$
is essentially self-dual if and only if $\lambda_i + \lambda_{n+1-i} =
\lambda_j + \lambda_{n+1-j}$ whenever $1\le i,j\le n$.
\medskip

\noindent
{\bf Definition 8.2.} We say that a point $k = (k_1,k_2,\ldots,k_n)\in
\cW(K)$
of the weight space is essentially self dual if $k_i + k_{n+1-i} = k_j +
k_{n+1-j}$ whenever $1 \le i,j\le n$. We let
$
\cW^{esd}
$
be the set of all essentially self-dual weights in $\cW$.
We also define
$$
\cXoe := \{\,x\in \cXoa\,|\,\wt(x) \in \cW^{esd}\}.
$$
\medskip

We note that $\cW^{esd}$ is a closed subspace of $\cW$ of dimension
$[n/2]+1$
and that $\cXoe$ is therefore a closed subspace of $\cXoa$ of dimension
$\le
[n/2]+1$.
\medskip

\noindent
{\bf Theorem 8.3.} {\it
Every cuspidal point of $\cXoa$ is contained in
$\cXoe$. }
\medskip

\noindent
{\bf Proof.} This is immediate from the Borel-Wallach theorem and the
fact that
$\cXoe$ is Zariski closed.
\medskip

We now suppose $n=3$.
We will use the notation $\lambda(g) = (2g,g,0)$ for $g \ge 0$.
Any essentially selfdual highest weight is of the form $\lambda(g)
\otimes \text{det}^j$ for some
$g \ge 0$ and $j \in \z$.  \medskip

\proclaim{Theorem 8.4} Let $g_0\ge 0$ and set $\lambda_0 =
\lambda(g_0)$. Let
$x_0\in\cXoa(\oQ_p)$ be a quasi-cuspidal arithmetic eigenpacket of
weight
$\lambda_0$ and suppose $x_0$ is not arithmetically rigid. Then for
every
integer
$g\equiv g_0$ modulo $p-1$ (or modulo $2$ if $p=2$), there is a cuspidal automorphic
point
$x\in\cXoa(\oQ_p)$ of weight
$\lambda(g)$ that is congruent to $x_0$.
\endproclaim

\demo{Proof} Since $x_0$ is not arithmetically rigid and let $Z$ be an
irreducible component of $\cXoa$ containing $x_0$. Then Corollary 4.8
and Theorem 4.9
tell us that $Z(\tilde x_0)$ is properly contained
in $Z$. In particular, we know that
$\dim Z \ge 2$. On the other hand, since $x_0$ is
quasi-cuspidal, every arithmetic point on $Z$ is quasi-cuspidal and is
therefore
cuspidal (theorem 6.4). So theorem 8.3 tells us $Z\subseteq \cXoe$. But
since
$n=3$, we know
$\cXoe$ has dimension $\le 2$. It follows that $Z$ is an irreducible
component of
$\cXoe$ and that $\dim Z = 2$. Since $\Roa$ is finite over $\Lambda$,
the image of
$Z$ under the weight map is a Zariski-closed irreducible
subset of $\cW^{esd}$ also having dimension $2$. But $\cW^{esd}$ is
$2$-dimensional
and therefore the image of $Z$ is precisely the connected component of
$\cW$
containing $\lambda_0$.

But if $g$ is a non-negative integer $g\equiv g_0$ modulo $p-1$ (or modulo $2$ if $p=2$), then
$\lambda(g)$
is in the same connected component of $\cW$ as $\lambda_0$. Thus, there
is a point
$x\in Z(\oQ_p)$ lying over $\lambda(g)$. By
lemma 4.3, $x$ is arithmetic. Clearly $x$ and $x_0$ are congruent (see
definition 4.5 and the accompanying remarks). Hence
$x$ is quasi-cuspidal and is therefore cuspidal automorphic of weight
$\lambda(g)$. This completes the proof of the theorem.
\enddemo

\medskip

The following corollary is immediate from the theorem and the
definitions.  Recall from the introduction that we denote by $V_g$ the representation $V_{\lambda(g)}$ of highest weight $\lambda(g)$.

\medskip

\proclaim{Corollary 8.5} Let $\varphi_0$ be a quasi-cuspidal arithmetic
eigenpacket of
weight
$\lambda(g_0)$ occurring in $H^3(\Gamma,V_{g_0}(K))$ and suppose
$\varphi_0$ is
not arithmetically rigid.
Then for every integer
$g\equiv g_0$ modulo $p-1$ (or modulo $2$ if $p=2$), there is a finite extension $L/K$ and a
cuspidal
automorphic eigenpacket $\varphi$ occurring in
$H^3(\Gamma,V_g(L))$
such that $\varphi$ is congruent to $\varphi_0$.
\endproclaim

\medpagebreak

\heading{Section 9:~~Examples of $p$-adically rigid classes via
computation}
\endheading

In this section, we assume $n=3$. According to Corollary 8.5, to check
that a given
quasi-cuspidal arithmetic eigenpacket is arithmetically rigid, it
suffices to show
$H^3_{cusp}(\Gamma,V_h(\C)) = 0$ for suitable $h$.

For any cohomology group $H^\ast()$ with coefficients in a module which
is a vector space
over a field $k$, let $h^\ast()$ denote its dimension over $k$. Let
$s_w(N)$ denote the
dimension of the space of classical holomorphic cuspforms of weight
$w$, level $N$ and trivial nebentype for
$GL(2)/\Q$.

\proclaim{Lemma 9.1} Let $q$ be prime. Then if $h$ is even,
$$
h^3(\Gamma_0(q),V_h(\C)) = h^3_{cusp}(\Gamma_0(q),V_h(\C)) +
2s_{h+2}(q) + 3
$$
and if $h$ is odd,
$$
h^3(\Gamma_0(q),V_h(\C)) = h^3_{cusp}(\Gamma_0(q),V_h(\C)).
$$
\endproclaim

\demo{Proof} This follows from the fact that $h^3_{cusp}(\Gamma,E) +
h^3(\partial Y,
E) = h^3(\Gamma,E)$ (see section 4). To compute the cohomology of the
Borel-Serre
boundary, take $\Gamma_0(q)$-invariants in
Theorem 7.6 (4) p. 116 in [Schwermer 83] (where $m = q$). Compare
Proposition 2.3 and Lemma 3.9 of [Ash-Grayson-Green 84].
\enddemo

Rather than compute over $\C$, we compute over $\z/r$ for some prime
$r$.   In fact we took $r = 149, 151$ or $163$. The following lemma
follows immediately
from the Universal Coefficient Theorem:

\proclaim{Lemma 9.2} Let $\Gamma = \Gamma_0(q)$ for some odd prime $q$,
and let $r$ be a
prime. Then
$h^3(\Gamma,V_h(\z/r)) = h^3(\Gamma,V_h(\C)) + e$ for some non-negative
integer $e$.
\endproclaim

\proclaim{Lemma 9.3} Let $\Gamma = \Gamma_0(q)$ for some odd prime $q$,
$h\ge 0$ be an
even integer, and
$r$ be a prime. If $h^3(\Gamma,V_h(\z/r)) = 2s_{h+2}(q) + 3$ then
$h^3_{cusp}(\Gamma,V_h(\C)) = 0$ and $e=0$.
\endproclaim

\demo{Proof} This follows from the previous two lemmas.
\enddemo

  From [Ash-Grayson-Green 84] we derive the following two examples of
$p$-adically rigid cohomology
classes.
(On the other hand, we have no glimmer of any non-rigid
non-essentially-selfdual
cohomology classes, consistent with the conjecture we made in the
introduction.)

\proclaim{Theorem 9.4}

(1) Let $\varphi$ be
either of the two eigenpackets in $H^3_{cusp}(\Gamma_0(89),\C).$ Then
$\varphi$ is $3$-adically arithmetically rigid.

(2) Let $\psi$ be
either of the eigenpackets in $H^3_{cusp}(\Gamma_0(61),\C).$ Then
$\psi$ is
$5$-adically arithmetically rigid up to twist.
\endproclaim

\demo{Proof} (1) A simple calculation from the Hecke data given in
[Ash-Grayson-Green 84]
shows that for $p=3$, $\varphi$ is 2-quasicuspidal and ordinary.
Therefore by Corollary 8.5, if $\varphi$ were not 3-adically
arithmetically rigid,
there would be nontrivial cuspidal cohomology in
$H^3_{cusp}(\Gamma_0(89),V_2(\C))$. But the computations
reported upon in the introduction to this paper say that
$H^3_{cusp}(\Gamma_0(89),V_2(\C)) = 0$.

(2) Similarly, for $p=5$, $\psi$ is 2-quasicuspidal and ordinary. Since
our computations say that $H^3_{cusp}(\Gamma_0(61),V_4(\C)) = 0$, the
same
reasoning shows that $\psi$ is
$5$-adically arithmetically rigid up to twist.
\enddemo

\heading{References}
\endheading

[Ash-Doud-Pollack 02] Ash, A., Doud, D. and Pollack, D., Galois
representations
with conjectural connections to arithmetic cohomology, Duke Math. J.
112 (2002) 521-579.

[Ash-Grayson-Green 84] Ash, A., Grayson, D. and Green, P., Computations
of cuspidal cohomology
of congruence subgroups of $SL(3,\z)$, J. Number Theory, 19 (1984)
412-436.

[Ash-Stevens 86] Ash, A. and Stevens, G., Cohomology of arithmetic
groups and congruences between systems
of Hecke eigenvalues, J. f. d. reine u. angew. Math., 365 (1986)
192-220.

[Ash-Stevens 97] Ash, A. and Stevens, G., $p$-adic deformations of
cohomology classes of
subgroups of $GL(n,\z)$, Collect. Math. 48 (1997) 1-30.

[Berstein-Zelevinsky 77] Bernstein, I. N. and Zelevinsky, A. V, Induced
representations of reductive {$p$}-adic groups. {I}, Ann. Sci. \'Ecole
Norm. Sup. (4), 10 (1977) 441-472

[Borel 83] Borel, A.,
Regularization theorems in Lie algebra cohomology. Applications.
Duke Math J. 50 (1983) 605-623.

[Borel-Wallach 00] Borel, A. and Wallach, N.,
Continuous cohomology, Discrete Subgroups, and Representations of
Reductive
Groups, 2nd ed., Math. Surv. and Mono. 67, American Math. Soc.,
Providence, 2000.

[Cartier 79] Cartier P.,
Representations of $\frak p$-adic groups: a survey, in Proc Symp Pure
Math
33 (1979) part 1, 111-155.

[Clozel 90] Clozel, L., Motifs et Formes Automorphes, in Automorphic
Forms, Shimura Varieties and L-functions, vol. 1, ed. L. Clozel and J.
S.
Milne, Academic Press, Boston 1990.

[Coleman-Mazur 98] Coleman R., Mazur, B., The Eigencurve, in
Galois representations in arithmetic algebraic geometry (Durham, 1996),
1--113,
London Math. Soc. Lecture Note Ser., 254, Cambridge Univ. Press,
Cambridge, 1998.

[Emerton 04] Emerton, M., On the interpolation of systems of eigenvalues
attached to automorphic Hecke eigenforms. Preprint, 2004.

[de Jong 96] de Jong, A. J., Crystalline Dieudonn\'e module theory via
formal and
rigid geometry. Inst. Hautes \'Etudes Sci. Publ. Math. No. 82 (1995),
5--96 (1996).

[Hida 94]  Hida, H., $p$-adic ordinary Hecke algebras for $GL(2)$, Annales de L'Institute
Fourier 44 (1994) 1289-1322.

[Kudla 94] Kudla, Stephen S., The local {L}anglands correspondence: the
non-{A}rchimedean case, in Proc Symp Pure Math 55 (1994) part 2,
365-391

[Schwermer 83] Schwermer,
J., Kohomologie airthmetisch definierter Gruppen und Eisensteinreihen,
Lecture Notes in
Mathematics 988, Springer Verlag, New York, 1983.

[Shimura 71] Shimura, G., Introduction to the arithmetic theory of automorphic
 functions,Publications of the Mathematical Society of Japan, No. 11.
 Iwanami Shoten, Publishers, Tokyo, 1971

\enddocument